\documentclass[10pt,reqno]{amsart} 
\usepackage{amsthm,amsmath,amsfonts,amssymb, hyperref,cases,color}
\usepackage{tikz}
\usepackage[capitalise]{cleveref}
\usepackage[T1]{fontenc}
\usepackage[utf8]{inputenc}
\usepackage{enumitem}

\everymath{\displaystyle}

\numberwithin{equation}{section}
\newtheorem{Theorem}{Theorem}[section]
\newtheorem{Lemma}[Theorem]{Lemma}
\newtheorem{Proposition}[Theorem]{Proposition}
\newtheorem{Corollary}[Theorem]{Corollary}
\newtheorem{Definition}[Theorem]{Definition}
\newtheorem{Remark}[Theorem]{Remark}

\renewcommand{\leq}{\leqslant}

\renewcommand{\geq}{\geqslant}

\DeclareMathOperator{\Span}{span}

\renewcommand{\Re}{\operatorname{Re}}

\def\C{{\mathbb{C}}}

\def\N{{\mathbb{N}}}

\def\R{{\mathbb{R}}}

\def\Z{{\mathbb{Z}}}

\title[Short-time estimates for a real moment problem]{Short-time estimates for a real moment problem with Two-Term Weyl Spectral Law}
\author{R\'emi Buffe}\author{Alessandro Duca}

\thanks{Universit\'e de Lorraine, CNRS, IECL, F-54000 Nancy, France.}

\begin{document}

\begin{abstract}
In this work, we first study the solvability of moment problems involving real exponentials $e^{\lambda_k t}$ and provide explicit estimates of the associated control cost. The result holds when the increasing sequence of distinct real numbers $(\lambda_k)_{k\in\N^*}$ satisfies a suitable two-term Weyl asymptotic law, without imposing any uniform spacing condition on blocks of its elements. We then deduce a corresponding controllability result for a linear control problem. Next, we present an exponential family fitting our hypotheses that cannot be treated by existing results of this type. Finally, we show how to deduce new exact controllability results for suitable fractional bilinear heat equations in higher-dimensional domains.\end{abstract}

\maketitle

\tableofcontents
\section{Introduction and main results}
\subsection{Moment problems}

Different problems arising in applied mathematics can be formulated in terms of the solvability of a so-called moment problem. Let $\mathsf H$ be a Hilbert space defined over a scalar field $\mathsf F$, and let $(\xi_k)_{k\in I}$ be a family of elements of $\mathsf H$, with $I\subseteq\Z$. The solvability of the associated moment problem consists in determining whether, for some specific given sequence of scalars $(x_k)_{k\in I}\subset \mathsf F$, there exists $h\in \mathsf H$ such that
$$\langle \xi_k, h\rangle_{\mathsf H}=x_k,\quad\forall k\in I.$$

For instance, let $\mathsf H=L^2(0,2\pi;\C)$ and $\xi_k= e^{i k t}$ with $k\in I=\Z$, the moment problem reads
$$\int^{2\pi}_0 e^{-i k t} h(t) dt =x_k,\quad \forall k\in \Z,$$
for $(x_k)_{k\in \Z}\in \ell^2(\Z,\C)$. In this case, the solution $h$ is uniquely reconstructed from the coefficients $(x_k)_{k\in \Z}$ by considering them as the Fourier coefficients of $h$ (up to a normalization factor). More generally, this type of construction applies whenever $(\xi_k)_{k\in I}$ is an orthogonal basis of $\mathsf H$.
\medskip

In this work, we consider the Hilbert space $\mathsf H=L^2(0,T;\R)$, with $T>0$, together with a family of real exponentials $$\xi_k(t)= e^{\lambda_k t},\quad k\in I=\N^*,\;t\in(0,T),$$ generated by an increasing sequence of distinct positive numbers $(\lambda_k)_{k\in \N^*}\subset\R$. We investigate the solvability of the following moment problem: given $x=(x_k)_{k\in \N^*}\subset\R$ in a suitable normed sequence space $\mathsf X$, we ask whether it is possible to find $h$ such that
\begin{equation}\label{mom}\int^{T}_0 e^{\lambda_k t} h(t) dt =x_k,\quad \forall k\in \N^*.\end{equation}
In addition, we estimate the cost associated with \eqref{mom}, defined as the optimal constant $C_T>0$ such that
$$\|h\|_{L^2(0,T)}\leq C_T \| x\|_{\mathsf X},$$
given $h$ solving \eqref{mom}. In this setting, the functions $\xi_k$ are not guaranteed to form either a basis or an orthonormal system, and the solvability of the moment problem requires techniques that differ from the standard orthogonal expansion framework.\\

For a real sequence of non-decreasing positive numbers $\lambda=(\lambda_k)_{k\in\N^*}$, we define its counting function
\[
\mathcal N_{\lambda}(\Gamma) = \#\{k\in\N^{*},\; \lambda_{k}\leq \Gamma\}.
\] 
Let $0<b<a$ and $p\in\R$. We consider a sequence $\mu:=(\mu_{k})_{k\in\N^{*}}\subset\R$ satisfying the following properties
\begin{enumerate}
\item The sequence
\begin{equation} \label{DefOfMu}
\mu \text{ is increasing and }1\leq \mu_1, 
\end{equation}
\item There exists $c_{W,1},c_{W,2}>0$ such that for any $k\in\N^{*}$ and any $\Gamma\geq\mu_1$,
\begin{equation}\label{WeylLaw}
c_{W,1} \Gamma^{a}-c_{W,2}\Gamma^{b} \leq \mathcal N_{\mu}(\Gamma) \leq c_{W,1} \Gamma^{a}+c_{W,2}\Gamma^{b},
\end{equation}
with the assumption
\begin{equation}\label{AssumptionOnAB}
\frac{1}{2} \leq a < \frac{5}{8}\text{ and }0< b<\frac{1}{2}.
\end{equation}
\item There exist $c_{w}>0$ and $p\in\R$ such that for any $k\in\N^{*}$
\begin{equation}\label{Weakgap}
\mu_{k+1}-\mu_{k} \geq \exp\left( - c_w \mu_{k}^{\frac{1}{2}} \right).
\end{equation}
\end{enumerate}
The requirement $a\geq \tfrac{1}{2}$ in \eqref{AssumptionOnAB} arises from the technique used in the main part of the article. However, this restriction is not intrinsic and can be relaxed to $0<a<\tfrac{1}{2}$ (see Appendix \ref{WeakerAccumulationSection}). Also note that the condition \eqref{DefOfMu} is not essential for the applications, and shall be used for convenience. In general, an appropriate change of variables allows one to shift the spectrum of the operator under consideration and then deal with spectra whose first eigenvalue is smaller than $1$. \\

From \eqref{WeylLaw}, it is immediate that there exists $K_{W}:=K_W(a,b,c_{W,1},c_{W,2})>0$ such that for any $k\in\N^{*}$ and any $\Gamma>0$
\begin{equation}\label{WeylLaw1Term}
\frac{1}{K_{W}}k^{\frac{1}{a}} \leq \mu_{k} \leq K_{W}k^{\frac{1}{a}},\quad \frac{1}{K_{W}^{a}}\Gamma^{a} \leq \mathcal N_{\mu}(\Gamma) \leq K_{W}^{a} \Gamma^{a}.
\end{equation}
For $0<b<a$ and $\delta\in(0,1)$, we define 
\begin{equation}\label{DefOfAlpha}
    \alpha_{\varepsilon,\delta}=\max\left( \frac{1-\varepsilon}{\varepsilon},\frac{1-\delta}{\delta} \right)
    \end{equation}
    where $\varepsilon =\tfrac{1}{4}\min \left(5-8a,1-2b\right)$.

  \subsubsection{On the solvability of a moment problem}
The main result of this article is as follows.
\begin{Theorem}\label{SolvabilityOfTheMoment}
 Let $\mu=(\mu_{k})_{k\in\N^{*}}$ be an increasing sequence of positive numbers satisfying \eqref{DefOfMu}, \eqref{WeylLaw}, \eqref{AssumptionOnAB} and \eqref{Weakgap}. Let $z=(z_k)_{k\in\N^*}\in \ell^2$ such that there exists $C_z>0$ such that for any $k\in\N^*$,
\begin{equation}\label{ConditionZk}
|z_k|\geq \exp \left( -C_z\mu_k^{1-\delta} \right)
\end{equation}
for some $\delta\in(0,1)$. Then, there exists $C:=C(a,b,\varepsilon,\delta,z,\mu,c_{W,1},c_{W,2},c_w)>0$ such that for any $T\in(0,1)$ and for any $x=(x_{k})_{k\in\N^{*}}\subset \R$ satisfying $zx \in\ell^2$, there exists $u\in L^{2}(0,T)$ such that
\[
x_{k}=\int_{0}^{T}e^{t\mu_{k}}u(t) dt,
\]
with the following estimate
\[
\|u\|_{L^{2}(0,T)} \leq C\exp\left( \frac{C}{T^{\alpha_{\varepsilon,\delta}}} \right) \| zx \|_{\ell^{2}}.
\]
\end{Theorem}

The proof is given in Section \ref{ProofOfTheMainResult}. An important example of a sequence that fulfills the hypotheses of Theorem \ref{SolvabilityOfTheMoment} but fails to satisfy any uniform block-spacing condition is given in Section \ref{ExampleSectionMoment}.
\begin{Remark} Few remarks are in order.
\begin{enumerate}
 \item The condition $a\geq 1/2$ seems somewhat artificial, as one might expect the situation to be more favorable whenever $0<b<a<\tfrac{1}{2}$. As our strategy is adapted to the case $\tfrac{1}{2}\leq a<\tfrac{5}{8}$, Theorem \ref{SolvabilityOfTheMoment} is stated for this range of parameter. In Appendix~\ref{WeakerAccumulationSection}, we provide the counterpart of Theorem~\ref{SolvabilityOfTheMoment} whenever $0<b<a<\tfrac{1}{2}$. However, owing to the method employed in the proof, the resulting control cost does not seem to be optimal.
 \item The condition $a<\tfrac{5}{8}$ appears as a barrier in the method we use here. However, in light of \cite[Theorem 2.79]{coron}, one may expect that these results might also hold for $\tfrac{5}{8}\leq a<1$. 
 \end{enumerate}
\end{Remark}

\subsubsection{Null-controllability of an abstract linear scalar control problem}\label{NullControlLinearSection}
    
Let $(\mu_k)_{\in \N^*}$ be a increasing sequence of positive numbers satisfying \eqref{DefOfMu}, \eqref{WeylLaw}, \eqref{AssumptionOnAB},\eqref{Weakgap}. Let $\mathcal H$ be a real separable Hilbert space endowed with a scalar product $\langle\cdot,\cdot\rangle_{\mathcal H}$, admitting $(\phi_k)_{k\in\N^*}$ as a Hilbert basis. We define the following operator
\begin{equation}\label{DefOfA}
\begin{array}{rcl}
\mathcal A : D(\mathcal A) \subset \mathcal H & \longrightarrow &\mathcal H \\
\sum_{k\in \N^*}\psi_k \phi_k& \longmapsto & \sum_{k\in \N^*}\psi_k \mu_k\phi_k,
\end{array}
\end{equation}
with domain 
\[D(\mathcal A)=\{\psi\in \mathcal H,\;\sum_{k\in\N^*} |\langle \psi,\phi_k\rangle_{\mathcal H}\mu_k|^2<+\infty\}.\]
Here, we assume that $\mathcal B\in \mathcal H$ satisfies: there exists $C_{\mathcal B}>0, \delta\in(0,1)$ such that for any $k\in\N^*$ 
\begin{equation}\label{AssumptionOnB}
    |\langle \mathcal B,\phi_k\rangle_{\mathcal H}|\geq \exp\left(-C_{\mathcal B} \mu_k^{1-\delta} \right).
\end{equation}
We are interested in the following control problem 
\begin{equation}\label{lin2}
\begin{cases}
\xi'(t)+\mathcal A\xi(t)+\mathcal Bv(t)=0 & t\in(0,T)\\
\xi(0)=\xi_0\in \mathcal H,
\end{cases}
\end{equation}
where $v\in L^2(0,T)$ is the control.
\begin{Theorem}\label{NullControlTheorem}
There exists $C:=C(a,b,\varepsilon,\delta,\mathcal B,\mu,c_{W,1},c_{W,2},c_w)>0$ such that for any $\xi_0\in \mathcal H$ and any $T\in(0,1)$, there exists $v\in L^2(0,T)$ such that the unique mild solution of \eqref{lin2} satisfies
\[
\xi(T)=0,
\]
    with control cost estimate
    \begin{equation}\label{NullControlCost}
    \| v \|_{L^2(0,T)} \leq C\exp \left(\frac{C}{T^{\alpha_{\varepsilon,\delta}}}\right) \|\xi_0\|_{\mathcal H},
    \end{equation}
    where $\alpha_{\varepsilon,\delta}$ is defined in \eqref{DefOfAlpha}.
\end{Theorem}
The proof is postponed in Section \ref{ProofOfTheSecondTheorem}.
\begin{Remark}
\begin{enumerate}
\item As for Theorem \ref{SolvabilityOfTheMoment}, we restrict our attention in Theorem \ref{NullControlTheorem} to the case $\tfrac{1}{2} \leq a < \tfrac{5}{8}$. Obviously, regarding the solvability of the moment problem, the controllability result can be extended to the case $0\leq a < \tfrac{1}{2}$. We provide the corresponding result in Appendix~\ref{WeakerAccumulationSection} for completeness, although the resulting control seems far from being optimal. 
 \item Conditions \eqref{ConditionZk} and \eqref{AssumptionOnB} are clearly related and are classical assumptions in the literature, imposed to guarantee null controllability in any positive time. Assuming $\delta = 0$ deteriorates the control operator and introduces a positive minimal control time (see \textit{e.g} \cite{ammar1,Dolecki}).
 \item Condition \eqref{DefOfMu} in Theorem \ref{NullControlTheorem} may be shifted to $\mu_1\geq -\sigma$, $\sigma\geq 0$, using the classical change of unknown $\widetilde \xi(t)= e^{(\sigma+\mu_1) t}\xi(t)$ and $\widetilde v (t) = e^{(\sigma+\mu_1) t}v(t).$
\end{enumerate}
\end{Remark}

\subsubsection{Literature overview}
The solvability of moment problems involving real exponentials to establish controllability properties for parabolic equations goes back to the seminal work of Fattorini and Russel \cite{fattorini,fattorini1}. There, spectral assumptions are made, \textit{e.g}
\begin{equation}\label{SommableSequence}
\sum_{k=1}^{+\infty} \frac{1}{\mu_k} <+\infty,
\end{equation}
and
\begin{equation}\label{StrongGap}
    \inf_{k\in\N^*}|\mu_{k+1}-\mu_k|\neq 0. 
\end{equation}

In higher dimensions, Fursikov and Imanuvilov \cite{FI} and Lebeau and Robbiano \cite{Lebeau_Robbiano} proved the null-controllability of the heat equation using Carleman estimates, with localized controls depending on both space and time. These results show that no minimal time is required, that is, the solution can be driven to rest in an arbitrarily small time. However, Dolecki \cite{Dolecki} identified a specific control operator for which a positive minimal time is necessary. This is now a very active field of research \cite{morgan,ABGBM, ammar1,BBGBO1,BBM2,Boyer,CMV1,CMV2,citGBO}, where authors have explored how spectral conditions and control operator properties lead to the existence or non-existence of a minimal time in general abstract settings. 

\medskip

In \cite{ammar1}, under the hypothesis that the eigenvectors forms a Riesz basis, the minimal time needed to have a null-control depends on the condensation index of the spectral sequence $\mu$ and the spectral properties of the control operator. In fact, they proved that if the control operator has sufficiently good properties (namely, if it allows to reach all frequencies with small control cost), then the minimal time is equal to the condensation index of the spectral sequence (see \cite{ammar1,Bernstein,Shackell} or Section \ref{RemarkCondensationIndexSection} below for a definition). In this setting, the authors do not provide an explicit dependence of the control cost with respect to small times. 

\medskip

In \cite{BBM2,Boyer}, the authors succeeded to solve in arbitrarily short time a block moment problem where the spectral sequence admits some condensation, along with condensation of the associated (generalized) eigenvectors. There, and in the recent article \cite{citGBO}, small time estimates on the control cost are provided, under the assumption of \eqref{SommableSequence}, and where \eqref{StrongGap} is replaced by the existence of uniformly bounded block ensuring the spectral gap, of the form 
\begin{equation}\label{BlockGapCondition}
\exists p\in\N^*,\rho>0 \text{  s.t. }\# ( \mu \cap [\lambda,\lambda+p])<\rho,\quad \forall \lambda>0.
\end{equation}
Finally, recent progress \cite{ABGBM} has been made toward the solvability of the moment problem arising from controlled equations in greater dimension than one, that is when \eqref{SommableSequence} is not satisfied. \\

Here, we are concerned with small time null-controls costs. Our approach addresses cases where \eqref{SommableSequence} holds but \eqref{BlockGapCondition} may fail (see Section \ref{ExampleSectionMoment}), \textit{i.e.},
\[
\{M\in\N^*,\; \inf_{k\in\N^*} |\mu_{k+M}-\mu_k|\neq 0 \}=\emptyset. \]
We initially follow the framework of \cite{citGBO} (see Proposition \ref{GBO} below). In this Proposition, we observe that the constants in the exponential estimates can be improved in the following manner: by allowing $\mathfrak{p}_2$ and $\mathfrak{ g p_1}^{-1} $ to depend slightly on the frequency, we can weaken the assumptions on the Weyl asymptotics and spectral gap. This is achieved through a Lebeau-Robbiano strategy: we first control the frequencies below a cut-off $\Lambda>0$ (where $\mathfrak{p}_2$ and $\mathfrak{ g p_1}^{-1}$ may grow slowly with $\Lambda$), and then leverage the dissipation of the heat equation to complete the construction of the final null-control.
\subsubsection{A remark on the condensation index} \label{RemarkCondensationIndexSection}
We introduce the following definition \cite{ammar1,Bernstein,Shackell}.
\begin{Definition}
Let $\lambda=(\lambda_k)_{k\in\N^*}\subset \C^*$ such that $\lambda^{-1}\in\ell^1$ $\lambda_i\neq\lambda_j$ for any $i,j\geq 1$, and such that there exists $\alpha>0$ such that $\Re(\lambda_k)\geq \alpha |\lambda_k|$ for any $k\geq1$. We define its condensation index
\[
c(\lambda) = \limsup_{k\to+\infty} \frac{\log\left( \frac{1}{E'(\lambda_k)}\right) }{\Re(\lambda_k)},
\]
where $E$ is the interpolating function
\[E(z)=\prod_{k\in\N^*}\left( 1 -\left | \frac{z}{\lambda_k}\right|^2\right),\qquad z\in\C.\]
\end{Definition}
We have proven the following result.
\begin{Corollary}
    Let $\mu\subset \R$ satisfying \eqref{DefOfMu}, \eqref{WeylLaw}, \eqref{AssumptionOnAB},\eqref{Weakgap}. Then
    \[
    c(\mu)=0.
    \]
\end{Corollary}
\begin{proof}
    As \eqref{lin2} is null-controllable for any $T\in(0,1)$, then the result follows as a consequence of \cite[Theorem 2.6]{ammar1}
\end{proof}

\subsubsection{A sequence with blocks of unbounded size} \label{ExampleSectionMoment}
We provide an example of a spectral sequence that fulfills the hypotheses of Theorem~\ref{SolvabilityOfTheMoment} yet fails to satisfy any uniform block spacing property. Let $\Omega=(0,1)^{2}$ be the unit square in two dimensions. Let $-\Delta_{D}$ be the associated Dirichlet Laplace operator and $\lambda=(\lambda_k)_{k\in\N^*}$ be its ordered nondecreasing sequence of eigenvalues. We consider
\[
\mu_k=\lambda_k^2- \frac{1}{ k}+2,\quad k\in\N^*.
\]
Let us show that $\mu_k$ satisfies \eqref{WeylLaw},\eqref{AssumptionOnAB},\eqref{Weakgap}. First, we know from the Weyl law on the square, that there exists $c_1>0$ such that for any $\Gamma>0,$
\[| \mathcal N_{\lambda^2}(\Gamma)-c_1\Gamma^{\frac{1}{2}} | = | \mathcal N_{\lambda}(\Gamma^{\frac{1}{2}})-c_1\Gamma^{\frac{1}{2}} | = \mathcal O (\Gamma^{\frac{1}{4}}).
\]
As a consequence,
\[| \mathcal N_{\mu}(\Gamma)-c_1\Gamma^{\frac{1}{2}} | = \mathcal O (\Gamma^{\frac{1}{4}}).
\]
Then, the sequence $\mu$ satisfies \eqref{WeylLaw} and \eqref{AssumptionOnAB} with $a=\frac{1}{2}$ and $b=\frac{1}{4}$. Furthermore, $\mu$ satisfies \eqref{Weakgap} since 
\[
\mu_{k+1}-\mu_k \geq  \frac{1}{k} \gtrsim \frac{1}{\mu_k^a}.
\]
Yet, from \cite{FSHH}, the sequence $\lambda$ enjoys the following property : for any $N\in\N^*$, there exists $k(N)\in\N^*$ such that 
\[
\#\{\ell\in \N^*,\; \lambda_\ell=\lambda_{k(N)}\} = N.
\]
This ensures $\{M\in\N^*,\; \inf_{k\in\N^*} |\mu^s_{k+M}-\mu^s_k|\neq 0 \}=\emptyset$ for any $s\in\R$.

\subsection{Bilinear controllability of parabolic equations}
\subsubsection{An abstract result}\label{intro_bili}
Let $(\mu_k)_{\in \N^*}$ be a increasing sequence of positive numbers satisfying \eqref{DefOfMu}, \eqref{WeylLaw}, \eqref{AssumptionOnAB},\eqref{Weakgap}, and let $\mathcal A$ be the operator given by \eqref{DefOfA} on a separate real Hilbert space $\mathcal H$.
We define the so-called $j-$th eigensolutions of $\mathcal A$ given by 
\begin{equation*}
\Phi_{j} (t)= e^{-\mu_{j}  t}\phi_{j} ,\quad t>0,\;\forall j\in\N^*.
\end{equation*}
Let $\mathfrak B\in \mathcal L(\mathcal H)$ be such that for any $j\in\N^*$, there exists $C_{\mathfrak B,j}>0, \delta\in(0,1)$ such that for any $k\in\N^*$ 
\begin{equation}\label{AssumptionOnFrakB}
    |\langle \mathfrak B\phi_j,\phi_k\rangle_{\mathcal H}|\geq \exp\left(-C_{\mathfrak B,j} \mu_k^{1-\delta} \right).
\end{equation}
We are interested in the following bilinear control problem 
\begin{equation}\label{main}
\begin{cases}
\psi'(t) +\mathcal A\psi(t)+v(t) \mathfrak B \psi(t) =0, &t\in(0,T),\\
\psi(0)=\psi_0\in \mathcal H,
\end{cases}
\end{equation}
where $v\in L^2(0,T)$ is the control. The well-posedness of \eqref{main} is a particular case of \cite[Theorem 2]{bms} and is recalled in Proposition \ref{well-posedness}. 
\begin{Definition}\label{ECTEDefinition}
 We say that \eqref{main} is small-time locally exactly controllable to eigensolutions if for any $T>0$ and any $j\in\N^*$, there exists $\delta>0$ such that for any $\psi_0\in L^2(\Omega,\R)$ with $$\|\psi_0-\phi_{j} \|_{L^2}<\delta,$$ there exists $v\in L^2(0,T)$ 
such that the unique mild solution to \eqref{main} satisfies 
$$\psi(T)=\Phi_{j}(T) .$$ 
\end{Definition}

We now state the application of the previous results to the local exact bilinear controllability to the eigensolutions.

\begin{Theorem} \label{BilinearControlTheo}
System \eqref{main} is small-time locally exactly controllable to eigensolutions.
\end{Theorem}

We provide in Section \ref{Section_Control_Simple_Spectrum} an explicit example based on the Bi-Laplace operator endowed with Dirichlet boundary conditions in dimension 2 where the hypothesis of Theorem \ref{BilinearControlTheo} are fulfilled (see Theorem \ref{th:local-rec}).

\subsubsection{Literature overview}

In the well-known result \cite{bms}, Ball, Marsden, and Slemrod established that multiplicative bilinear controls can introduce obstructions to the exact controllability of a system dynamics. Roughly speaking, when the control operator $\mathfrak B$ is bounded on $\mathcal H$, the authors showed that the reachable set is contained in a countable union of compact sets. This phenomenon is very general: it occurs not only to the bilinear heat equation, but also to the bilinear Schrödinger and wave equations. The negative result established in \cite{bms} is the primary motivation for studying reachable targets in evolution equations with bilinear controls. One can cite the important works \cite{Beauchard1,laurent} on the exact controllability of bilinear wave or Schr\"odinger equations, in the case where the control operator is not bounded.

\smallskip

Small-time local exact controllability to eigensolutions properties for the heat-like equation \eqref{main} was firstly derived by Alabau-Boussouira, Cannarsa and Urbani in \cite{acu1,acue} assuming the gap condition 
\begin{equation}\label{gap-intro}
\exists \gamma>0 ,\; \forall k\in\N^*,\quad |\sqrt{\lambda_{k+1}-\lambda_j}-\sqrt{\lambda_k-\lambda_j}|\geq \gamma>0
\end{equation}
and the spreading assumption
\begin{equation}\label{ipB-intro}
\langle \mathcal B\phi_j,\phi_j\rangle_{\mathcal H}\neq 0,\quad|\lambda_k-\lambda_j|^q|\langle \mathcal B\phi_j,\phi_k\rangle_{\mathcal H}|\geq b,\quad\forall\,k\neq j,
\end{equation}
with $p,b>0$. A direct application of their work is the controllability of the heat equation \eqref{main} in the one-dimensional case $\Omega=(0,1)$. For related problems and further developments on this topic, see also \cite{acu, cu}.

\smallskip

 These results have been extended in \cite{cdu} by the second author with Cannarsa and Urbani, where they proved the same local controllability to the eigensolutions for the equation \eqref{main} on network-type domains. In this framework, the spectrum of the Laplacian $-\Delta$ has different properties with respect to the unit interval as we only know the existence of $N\in N^*$ such that
\begin{equation}\label{gap_weak1}
\sqrt{\lambda_{k+N}}-\sqrt{\lambda_k}\geq\gamma>0.
\end{equation}

\smallskip 

In this work, as an application of Theorem \ref{SolvabilityOfTheMoment}, we are able to generalize these results to situations where condition \eqref{gap_weak1} is not verified (see the example in Section \ref{ExampleSectionMoment}). We then ensure the local controllability to the eigensolutions stated in Theorem \ref{BilinearControlTheo} by classical techniques employed in \cite{acu} (see also \cite{BeauchardMarbach,LTT}).

\subsection*{Acknowledgments} The authors were partially supported by the French National Research Agency (ANR), Project QuBiCCS (ANR-24-
CE40-3008). We thank Assia Benabdallah and Morgan Morancey for their valuable comments, which enhanced the clarity of the paper. We also thank the anonymous referees for the evaluation of this work.

\section{Solvability of a moment problem}

\subsection{Existence of biorthogonal families}

We shall use the following result, which is a particular case of \cite[Theorem 1.2]{citGBO}.
\begin{Proposition} \label{GBO}Let $(v_{k})_{k\in\N^{*}}$ be a sequence of increasing positive numbers such that there exists $\mathfrak p_{1},\mathfrak p_{2}\geq \mathfrak p_{0}>0$, $\alpha>0$ such that 
\begin{equation}\label{WeylGBO}
-\alpha + \mathfrak p_{1} \Gamma^{\frac{1}{2}} \leq \mathcal N_v(\Gamma) \leq \alpha+ \mathfrak p_{2}\Gamma^{\frac{1}{2}},\quad \Gamma>0.
\end{equation}
Assume that there exists $\mathfrak n\in\N^{*}$, $\mathfrak g>0$ and a non-decreasing sequence $(a_{k})_{k\in\N^{*}}$ of positive real numbers such that for any $k\in\N^{*}$,
\begin{equation}\label{GapsGBO}
\begin{cases}
v_{k+\mathfrak n}-v_{k}\geq \mathfrak g (2k\mathfrak n + \mathfrak n^{2}) \\
v_{k+1}-v_{k}\geq a_{k}^{-1}.
\end{cases}
\end{equation}
Then for any $T>0$ there exists $(\sigma_{k})_{k\in\N^{*}}\subset L^{2}(0,T)$ such that 
\[
\int_{0}^{T}\sigma_{k}(t)e^{v_{j}t}dt = \delta_{j,k},\quad \forall j,k\in\N^{*},
\]
and there exists $C>0$, depending only on $v_{1},\mathfrak p_{0},\alpha$, such that 
\begin{multline}\label{UpperboundForBiortho}
\|\sigma_{k}\|_{L^{2}(0,T)} \leq C  a_{k}^{\mathfrak n} \left( \frac{16+4\mathfrak g \mathfrak p_{2}^{2}}{\mathfrak g^{2} \mathfrak p_{1}^{4}}\right)^{2\mathfrak n -2}\\
\times \exp\Big( C \left(1+ \mathfrak n + T^{\frac{1}{2}}+ 16(\mathfrak p_{1}\mathfrak g )^{-2}+\mathfrak p_{2} \right)v_{k}^{\frac{1}{2}}+C\frac{(1+\mathfrak p_{2})^{2}}{T}-v_{k}T\Big).
\end{multline}
\end{Proposition}

\begin{proof}The result is a direct consequence of \cite[Theorem 1.2]{citGBO}. Indeed, (H2), (H3), and (H4) in \cite[Definition 1.1]{citGBO} are fulfilled. Next, conditions \eqref{GapsGBO} imply (H1). Let us show that (H5) is satisfied. First, we prove by induction that for any $j\in\N^*$
\begin{align*}\label{id-int}
(\mathcal P_j)\qquad v_{k+j\mathfrak n}-v_{k}\geq \mathfrak g ((k+j\mathfrak n)^2-k^2) ,\quad \forall k \in \N^*.
\end{align*}
Indeed, $(\mathcal P_1)$ holds from assumption \eqref{GapsGBO}. Moreover, assume that $(\mathcal P_j)$ holds for some $j\in\N^*$. Then, using $(\mathcal P_1)$ and $(\mathcal P_j)$, for any $k\in\N^*$,
\begin{align*}
    v_{k+(j+1) \mathfrak n}-v_{k}  =     v_{(k+j\mathfrak n)+\mathfrak n}- v_{k+j\mathfrak n} +v_{k+j\mathfrak n} - v_{k}  
     \geq  \mathfrak g ((k+(j+1)\mathfrak n)^2 - k^2), 
\end{align*}
and then the assertion holds for $j+1$, ensuring $(\mathcal P_{j+1})$. Hence $(\mathcal P_j)$ holds for any $j\in\N^*.$ Now, let $j,n\in\N^*$ be such that $j\mathfrak n < n < (j+1)\mathfrak n$. Then, for any $k\in\N^*$,
\begin{align*}
    \mu_{k+n}-\mu_k &\geq \mu_{k+j\mathfrak n}-\mu_k  
    \geq \mathfrak g (2kj\mathfrak n + j^2\mathfrak n^2) \\
   & \geq \frac{\mathfrak g}{4} (2k(j+1)\mathfrak n + (j+1)^2\mathfrak n^2)   \geq \frac{\mathfrak g}{4} (2kn+n^2).
\end{align*}
Thus, we have proven
\begin{align*}
v_{k+ n}-v_{k}\geq \frac{\mathfrak g}{4} ((k+ n)^2-k^2) ,\quad \forall k,n \in \N^*,\; n\geq \mathfrak n,
\end{align*}
which is (H5). Then, the sequence $v$ belongs to the class $\mathcal L(0,\tfrac{\mathfrak g}{4},\mathfrak n,\mathfrak p_0,\mathfrak p_1,\mathfrak p_2,\alpha)$, defined in \cite[Definition 1.1]{citGBO}. Hence, \cite[Theorem 1.2]{citGBO} ensures the existence of $(\widetilde\sigma_k)_{k\in\N^*}$, a bi-orthogonal family to the exponentials $(e^{-\nu_k t})_{k\in\N^*}$ in $L^2(0,T)$, that is
$$\delta_{k,j}=\int_0^T \widetilde\sigma_k(t)e^{-v_j t}dt.
$$
The bi-orthogonal family to the exponentials $(e^{\nu_k t})_{k\in\N^*}$ is simply obtained by $\sigma_k(t)=\widetilde \sigma_k(T-t)e^{-\nu_k T}$. Finally, the estimate \eqref{UpperboundForBiortho} is a direct consequence of the bound in \cite[Theorem 1.2]{citGBO}, adapted to our setting.

\end{proof}

\subsection{Gaps for spectral sequences}
In this section, we consider a sequence $\mu=(\mu_k)_{k\in\N^*}$ satisfying \eqref{DefOfMu}, \eqref{WeylLaw}, \eqref{AssumptionOnAB} and \eqref{Weakgap}. 
\subsubsection{An auxiliary spectral sequence}
Let $\Lambda\geq \mu_{1}$, and let $\kappa = K_{W}^{a}$. We introduce the following auxilliary sequence
\begin{equation}\label{DefOfNu}
\nu_{k}=\begin{cases}\mu_{k} & \text{ if }\mu_{k}\leq \Lambda \\ \kappa^{2}\Lambda^{1-2a}k^{2} & \text {if } \mu_{k}>\Lambda. 
\end{cases}
\end{equation}
For any $\Lambda\geq \mu_{1}$, we define $k^{*}_{\Lambda}$ as the unique integer such that 
\begin{equation}\label{DefOfkStarLambda}
\begin{cases} \mu_{k_{\Lambda}^{*}} \leq \Lambda \\ \mu_{k^{*}_{\Lambda}+1}>\Lambda.\end{cases}
\end{equation}
Here, we define some constants appearing in the spectral gaps as follows 
\begin{equation}\label{ValueOfConstants}
\begin{cases}
\varepsilon =\min \left( \tfrac{5-8a}{4},\tfrac{1-2b}{4} \right)\in (0,1) \\
\Lambda_{0} =   
\mu_{1}+(8K_{W}^{a}c_{W,1}c_{W,2})^{\frac{4}{1-2b}}+ (2K_{W}\mu_{k_{*}})^{4a}\geq 1\\
\theta  =\min \Big( \tfrac{1}{4(c_{W,1}+c_{W,2})^2} (2K_W^a+K_{W}^{2a})^{\frac{1-2a}{2}}  , (8k_{*}K_{W})^{-1} ,K_W^{-2a}\Big) >0,
\end{cases}
\end{equation}
where 
\begin{equation}\label{DefOfKStar}
    k_{*}= \left\lfloor \left( \tfrac{c_{W,2}K_W^{a+b}}{c_{W,1}}\right)^{\frac{a}{a-b}}\right\rfloor +1.
    \end{equation}
 Let $\Lambda>0$, and let $0<b<a$ satisfy \eqref{AssumptionOnAB}. We set
\begin{equation}\label{DefsOfGapConstants}
n_{\Lambda}=\lfloor \Lambda^{\frac{1-2\varepsilon}{2}}\rfloor+1,\quad \gamma_{\Lambda}=\theta \Lambda^{1-2a}.
\end{equation}

The main result of this section is as follows.

\begin{Proposition}\label{BiorthoAuxilliary} Let $\Lambda_{0},\varepsilon,\theta>0$ be defined in \eqref{ValueOfConstants}. Let $\Lambda\geq \Lambda_{0}$ and  $(\nu_k)_{k\in\N^{*}}$ be the sequence defined in \eqref{DefOfNu}, under the assumptions \eqref{DefOfMu}, \eqref{WeylLaw}, \eqref{AssumptionOnAB}, \eqref{Weakgap}. For any $T\in(0,1)$, there exists $(\sigma_{k})_{k\in\N^{*}}\subset L^{2}(0,T)$ such that
\[
\int_{0}^{T}\sigma_{k}(t)e^{\nu_{j}t}dt = \delta_{j,k},\quad j,k\in\N^{*},
\]
and there exists $C:=C(a,b,\varepsilon,\mu,c_{W,1},c_{W,2},c_w)>0$, independent on $\Lambda\geq \Lambda_{0}$ and $T\in(0,1)$, such that 
\begin{equation*}
\|\sigma_{k}\|_{L^{2}(0,T)}  \leq C \exp\left( C\big(\nu_{k}^{\frac{1}{2}}\Lambda^{\frac{1}{2}-\varepsilon}+\Lambda^{1-\varepsilon}+\frac{1}{T^{\alpha_\varepsilon}}\big)\right),\quad k\in\N^{*},
\end{equation*}
where $\alpha_\varepsilon= \frac{1-\varepsilon}{\varepsilon}$.
\end{Proposition}
Before proceeding to the proof (which is postponed at Section \ref{ProofOfProposition22}), we need to collect intermediate results. We recall \eqref{DefOfkStarLambda}.
\begin{Lemma}\label{BelowTheCut}

For any $\Lambda\geq \mu_{1}$ and any $k>k^{*}_{\Lambda}$, we have
\[
\nu_{k}>\Lambda.
\]
As a consequence, for any $\mu_{1}\leq \Gamma \leq \Lambda$, 
\[
\{k\in\N^{*},\; \mu_{k}\leq \Gamma\} = \{k\in\N^{*},\; \nu_{k}\leq \Gamma\}.
\]
\end{Lemma}
\begin{proof}
Let $\Lambda\geq \mu_{1}$ and $k>k_{\Lambda}^{*}$. Recalling \eqref{WeylLaw1Term}, we have 

\begin{equation*}
 \mu_{k} \leq K_{W}k^{\frac{1}{a}}.
\end{equation*}
As $\Lambda < \mu_{k_{\Lambda}^{*}+1}$, we deduce that 
\[
K_{W}^{-2a} \Lambda^{2a} < (k_{\Lambda}^{*}+1)^{2}.
\]
Hence,
\[
\nu_{k^{*}_{\Lambda}+1} = \kappa^{2} \Lambda^{1-2a} (k_{\Lambda}^{*}+1)^{2} > \kappa^{2} K_{W}^{-2a} \Lambda = \Lambda.
\]
Let $\mu_{1}\leq \Gamma\leq \Lambda$ and let $k\in \N^{*}$ be such that $\mu_{k}\leq \Gamma$. Then by definition \eqref{DefOfNu}, $\mu_{k}=\nu_{k}$, and then $\{k\in\N^{*},\; \mu_{k}\leq \Gamma\} \subset \{k\in\N^{*},\; \nu_{k}\leq \Gamma\}$. Conversely, let $k\in \N^{*}$ be such that $\nu_{k}\leq \Gamma$. Then $k\leq k_{\Lambda}^{*}$, and thus $\mu_{k}=\nu_{k}$.
\end{proof}

\begin{Lemma}\label{BoundForMu}
 Assume $\frac{1}{2}\leq a$. For any $\Lambda\geq \mu_{1}$, we have
\[
\mu_{k}\leq \kappa^{2} \Lambda^{1-2a}k^{2}=\nu_k,\quad \forall k\in\N^{*},\;\mu_{k}> \Lambda.
\]
\end{Lemma}
\begin{proof}
Let $\mu_{1}\leq \Lambda$. Using \eqref{WeylLaw1Term} ,
\[
 \mu_{k} =\kappa^{2}\mu_{k}^{1-2a} K_{W}^{-2a}\mu_{k}^{2a} \leq  \kappa^{2}\mu_{k}^{1-2a}k^{2} .
\]
As $\mu_{k}^{1-2a}\leq \Lambda^{1-2a}$, we obtain the result.

\end{proof}

\begin{Lemma}\label{BoundsOnCounting}
There exists $c_{1},c_{2}>0$ such that for any $\Lambda,\Gamma\geq \mu_{1}$, 
\[
c_{1}\Gamma^{\frac{1}{2}} \leq \mathcal N_{\nu}(\Gamma) \leq c_{2} \Gamma^{\frac{1}{2}}\Lambda^{\frac{2a-1}{2}} 
\]
\end{Lemma}

\begin{Remark}
Inspecting the proof below yields, for instance
\[
c_{1}=\frac{1}{2}\min \Big(  K_W^{-2a},K_W^{-a}\Big),
\quad
c_{2}=2\max\Big(K_W^a,K_W^{-a}\Big).
\]
In particular, we emphasize that the constants $c_{1},c_{2}$ do not depend on the relevant parameters.
\end{Remark}
\begin{proof}
\textbf{Step 1 :} Assume $\mu_{1}\leq \Lambda<\Gamma$. According to Lemma \eqref{BelowTheCut}, we have the following identities
\begin{align}
\mathcal N_{\nu}(\Gamma) & = \mathcal N_{\nu}(\Lambda) +\# \{ k\in\N^{*},\; \Lambda < \nu_{k}\leq \Gamma \}\nonumber \\ 
  & = \mathcal N_{\mu}(\Lambda) +\lfloor \kappa^{-1}\Lambda^{\frac{2a-1}{2}}\Gamma^{\frac{1}{2}}\rfloor  - \lfloor \kappa^{-1}\Lambda^{a}\rfloor \label{EqualityOnCountings}. 
\end{align}
Using \eqref{WeylLaw1Term}, there holds
\[
\mathcal N_{\nu}(\Gamma ) \leq K_{W}^{a}\Lambda^{a}+\kappa^{-1}\Lambda^{\frac{2a-1}{2}}\Gamma^{\frac{1}{2}} \leq 2\max(\kappa,\kappa^{-1})\Lambda^{\frac{2a-1}{2}}\Gamma^{\frac{1}{2}},
\]
and the upper bound follows. We now prove the lower bounds. Let $\Gamma^{*}= 4\kappa^{2}$.\\

\noindent
\textbf{Case 1 :} If $\mu_{1}\leq \Gamma\leq \Gamma^{*}$. Then 
\[
\mathcal N_{\nu}(\Gamma) \geq \mathcal N_{\mu}(\Lambda) \geq  K_{W}^{-a}\Lambda^{a} \geq  K_{W}^{-a} 
(\Gamma^{*})^{-\frac{1}{2}}\Gamma^{\frac{1}{2}}= \frac{1}{2\kappa^2}\Gamma^{\frac{1}{2}}. 
\]

\noindent
\textbf{Case 2 :} If $\Gamma^{*} < \Gamma$.  
 Using \eqref{EqualityOnCountings} and \eqref{WeylLaw1Term},
\[
\mathcal N_{\nu}(\Gamma ) \geq K_{W}^{-a}\Lambda^{a}- \lfloor K_{W}^{-a}\Lambda^{a} \rfloor+\lfloor\kappa^{-1}\Lambda^{\frac{2a-1}{2}}\Gamma^{\frac{1}{2}}\rfloor \geq \lfloor\kappa^{-1}\Lambda^{\frac{2a-1}{2}}\Gamma^{\frac{1}{2}}\rfloor.
\]
As $1\leq \mu_{1}\leq\Lambda$ 
then 
\[
\kappa^{-1}\Lambda^{\frac{2a-1}{2}}\Gamma^{\frac{1}{2}} \geq \kappa^{-1} (\Gamma^{*})^{\frac{1}{2}}\geq 2.
\]
This implies
\[
\lfloor\kappa^{-1}\Lambda^{\frac{2a-1}{2}}\Gamma^{\frac{1}{2}}  \rfloor \geq \frac{1}{2}\kappa^{-1}\Lambda^{\frac{2a-1}{2}}\Gamma^{\frac{1}{2}}\geq  \frac{1}{2\kappa}\Gamma^{\frac{1}{2}}.
\]
Summing up, we proved the result in the case where $\Gamma>\Lambda$.\\

\noindent  
\textbf{Step 2 :} Assume $\mu_{1}\leq \Gamma\leq\Lambda$. In that case, using Lemma \ref{BelowTheCut}, we have 
\[
\mathcal N_{\nu}(\Gamma)=\mathcal N_{\mu}(\Gamma), 
\]
and then from \eqref{WeylLaw1Term},
\[
\frac{1}{K_{W}^{a}}\Gamma^{\frac{1}{2}} \geq \frac{1}{K_{W}^{a}}\Gamma^{a} \leq \mathcal N_{\nu}(\Gamma) \leq K_{W}^{a} \Gamma^{a} \leq K_{W}^{a} \Lambda^{\frac{2a-1}{2}}\Gamma^{\frac{1}{2}},
\]
and the result follows.
\end{proof}

\subsubsection{Estimates for the block gap condition}
Recall \eqref{ValueOfConstants}. The purpose of this section is to prove the following result.

\begin{Lemma}\label{GapBlockLemma}
For any $\Lambda\geq \Lambda_{0}$, and any $k\in\N^{*}$,
\[
\nu_{k+n_{\Lambda}}-\nu_{k}\geq \gamma_{\Lambda}(2kn_{\Lambda}+n_{\Lambda}^{2}).
\]
\end{Lemma}
We begin with an immediate property to which we shall refer.
\begin{Lemma}\label{EstimatesOnConstants} We have
\begin{equation}
\gamma^{-2}_{\Lambda}\leq \theta^{-2}\Lambda^{\frac{1}{2}-\varepsilon},\quad  1\leq \Lambda.
\end{equation}
\end{Lemma}
Note that we used the condition $a<\tfrac{5}{8}$ here.
\begin{Lemma}\label{GapMu}
For any $\Lambda\geq \Lambda_{0}$,
\[
\mu_{k+n_{\Lambda}}-\mu_{k}\geq \gamma_{\Lambda}(2kn_{\Lambda}+n_{\Lambda}^{2}),\qquad \forall k\in\N^{*},\; 1\leq \mu_{k} \leq \Lambda.
\]
\end{Lemma}
\begin{proof}
From \eqref{WeylLaw}, we deduce for any $k\in\N^*$,
\begin{equation}\label{CrucialBound}
c_{W,1} \mu_{k}^{a}-c_{W,2}\mu_{k}^{b} \leq k \leq c_{W,1} \mu_{k}^{a}+c_{W,2}\mu_{k}^{b}.
\end{equation}
Let $\Lambda\geq \mu_1\geq 1$ and recall that $k_*$ is defined in \eqref{DefOfKStar}. We shall separate three cases.\\

\noindent
\textbf{Case 1.} If $k\in \{j\in\N^{*},\; 1\leq \mu_{j} \leq \Lambda\}\cap \{j\in\N^{*},\; j\geq k_{*}\}$. From \eqref{DefOfKStar}, the left hand side of \eqref{CrucialBound} is non-negative. Hence, Estimate \eqref{CrucialBound} yields
\[
c_{W,1}^{2}\mu_{k}^{2a}+c_{W,2}^{2}\mu_{k}^{2b}-2c_{W,1}c_{W,2}\mu_{k}^{a+b} \leq k^{2} \leq c_{W,1}^{2}\mu_{k}^{2a}+c_{W,2}^{2}\mu_{k}^{2b}+2c_{W,1}c_{W,2}\mu_{k}^{a+b}.
\]
Then, for any $n\in\N^{*}$,
\begin{multline*}
(k+n)^{2}-k^{2}\leq c_{W,1}^{2}(\mu_{k+n}^{2a}-\mu_{k}^{2a})\\ +c_{W,2}^{2}(\mu_{k+n}^{2b}-\mu_{k}^{2b})+2c_{W,1}c_{W,2}(\mu_{k+n}^{a+b}-\mu_{k}^{a+b}) +4c_{W,1}c_{W,2}\mu_{k}^{a}\Lambda^{b} .
\end{multline*}
On the other hand, due to \eqref{WeylLaw1Term},
\[
\mu^{a}_{k} \leq K^{a}_{W}k = K^{a}_{W} \frac{k}{2kn+n^{2}}((k+n)^{2}-k^{2}) \leq \frac{K^{a}_{W}}{n} ((k+n)^{2}-k^{2}) .
\]
In particular, with $n=n_{\Lambda}$ (defined in \eqref{DefsOfGapConstants}) and as $0<\varepsilon\leq\frac{1-2b}{4}$,
\begin{align*}
4c_{W,1}c_{W,2}\mu_{k}^{a}\Lambda^{b}  &\leq 4K^{a}_{W}c_{W,1}c_{W,2} \Lambda^{\frac{2b+2\varepsilon-1}{2}} ((k+n_{\Lambda})^{2}-k^{2}).
\end{align*}
We set $\Lambda_{1}=1+(8K^{a}_{W} c_{W,1}c_{W,2})^{\frac{4}{1-2b}}$, and then for any $\Lambda\geq\Lambda_{1}$,
\[
4c_{W,1}c_{W,2}\mu_{k}^{a}\Lambda^{b}  \leq \frac{1}{2} ((k+n_{\Lambda})^{2}-k^{2}).
\]
Hence we have proven, for any $\Lambda\geq \Lambda_{1}$, for any $k\in \{k\in\N^{*},\; 1\leq \mu_{k} \leq \Lambda\}\cap \{k\in\N^{*},\; k\geq k_{*}\}$,
\begin{multline*}
(k+n_{\Lambda})^{2}-k^{2} \leq 2c_{W,1}^{2}(\mu_{k+n_{\Lambda}}^{2a}-\mu_{k}^{2a})\\ +2c_{W,2}^{2}(\mu_{k+n_{\Lambda}}^{2b}-\mu_{k}^{2b})+4c_{W,1}c_{W,2}(\mu_{k+n_{\Lambda}}^{a+b}-\mu_{k}^{a+b}).
\end{multline*}
\begin{Lemma}
For $A>B\geq0$ and $y\geq 1$, we have
\[
x^{A}-y^{A} \geq x^{B}-y^{B},\quad \forall x\geq y.
\]
\end{Lemma}

\noindent
We use that lemma with $A=2a$, $B=2b$, $x=\mu_{k+n_{\Lambda}}$, $y=\mu_{k}$, to obtain
\[
2c_{W,2}^{2}(\mu_{k+n_{\Lambda}}^{2b}-\mu_{k}^{2b}) \leq 2 c_{W,2}^{2}(\mu_{k+n_{\Lambda}}^{2a}-\mu_{k}^{2a}).
\]
We also use that lemma with $A=2a$, $B=a+b$, $x=\mu_{k+n_{\Lambda}}$, $y=\mu_{k}$, to obtain
\[
4c_{W,1}c_{W,2}(\mu_{k+n_{\Lambda}-\mu_{k}}^{a+b}-\mu_{k}^{a+b}) \leq 4 c_{W,1}c_{W,2}(\mu_{k+n_{\Lambda}}^{2a}-\mu_{k}^{2a}),
\]
which yields
\[
(k+n_{\Lambda})^{2}-k^{2} \leq 2(c_{W,1}+c_{W,2})^{2}(\mu_{k+n_{\Lambda}}^{a}-\mu_{k}^{a})(\mu_{k+n_{\Lambda}}^{a}+\mu_{k}^{a}).
\]
\begin{Lemma}
Let $m>1$. Then 
\[
x-y\leq x^{1-m}(x^{m}-y^{m}),\quad \forall\; 0<y<x.
\]
\end{Lemma}

\noindent
Applying this Lemma with $m=a^{-1}$ (recall that $a^{-1}>\frac{8}{5}$ from \eqref{AssumptionOnAB}), $x=\mu_{k+n_{\Lambda}}^{a}$, $y=\mu_{k}^{a}$, we obtain
\begin{align*}
(k+n_{\Lambda})^{2}-k^{2}
& \leq 4(c_{W,1}+c_{W,2})^{2}\mu_{k+n_{\Lambda}}^{2a-1}(\mu_{k+n_{\Lambda}}-\mu_{k}).
\end{align*}
Yet, from \eqref{WeylLaw1Term}
\begin{align*}
(k+n_{\Lambda})^{2}-k^{2}  \leq 4K_{W}^{2a-1}(c_{W,1}+c_{W,2})^{2}(k+n_{\Lambda})^{\frac{2a-1}{a}}(\mu_{k+n_{\Lambda}}-\mu_{k}) 
\end{align*}
As $\frac{1}{2}\leq a$, we have
\begin{align*}
(k+n_{\Lambda})^{2}-k^{2} & \leq 4K_{W}^{2a-1}(c_{W,1}+c_{W,2})^{2}(K_{W}^{a}\Lambda^{a}+ \Lambda^{\frac{1}{2}-\varepsilon} +1)^{\frac{2a-1}{a}}(\mu_{k+n_{\Lambda}}-\mu_{k})  \\
& \leq 4K_{W}^{2a-1}(c_{W,1}+c_{W,2})^{2}(2+K_{W}^{a})^{\frac{2a-1}{a}} \Lambda^{2a-1}(\mu_{k+n_{\Lambda}}-\mu_{k}) 
\end{align*}
Hence, we have proven in Case 1 that for any $\Lambda\geq \Lambda_{0}$,
\begin{align*}
\Lambda^{1-2a}((k+n_{\Lambda})^{2}-k^{2}) \leq   4K_{W}^{2a-1}(c_{W,1}+c_{W,2})^{2}(2+K_{W}^{a})^{\frac{2a-1}{a}} (\mu_{k+n_{\Lambda}}-\mu_{k}),
\end{align*}
which yields the announced inequality.\\

\noindent
\textbf{Case 2.}  If $k\in \{k\in\N^{*},\; 1\leq \mu_{k} \leq \Lambda\}\cap \{k\in\N^{*},\; k< k_{*}\}$. On the one hand, using \eqref{WeylLaw1Term},
\[
\mu_{k+n_{\Lambda}}-\mu_{k}\geq \mu_{n_{\Lambda}}-\mu_{k_{*}} \geq \frac{1}{K_{W}}n_{\Lambda}^{\frac{1}{a}}-\mu_{k_{*}} \geq \frac{1}{2K_{W}}\Lambda^{\frac{1-2\varepsilon}{2a}},
\]
 for $\Lambda\geq 1+ (2K_{W}\mu_{k_{*}} )^{4a}$. On the other hand, $\forall \Lambda\geq 1$, 
\[
\gamma_{\Lambda}(2kn_{\Lambda}+n_{\Lambda}^{2}) \leq \theta\Lambda^{2-2a-2\varepsilon} 4k_{*}.
\]
Summing up, for any $\Lambda\geq  1+ (2K_{W}\mu_{k_{*}} )^{4a}$,
\[
\gamma_{\Lambda}(2kn_{\Lambda}+n_{\Lambda}^{2}) \leq \theta\Lambda^{2-2a-2\varepsilon} 4k_{*} \leq K_{W} 8k_{*}  \theta\Lambda^{\frac{1}{2a}(4a-4a^{2}-1-4a\varepsilon+2\varepsilon)} (\mu_{k+n_{\Lambda}}-\mu_{k}).
\]
As $4a-4a^{2}-1-4a\varepsilon+2\varepsilon= -4(a-\frac{1}{2})(a-\frac{1}{2}+\varepsilon)$ and $a\geq   \frac{1}{2},$ we deduce, that this quantity is nonpositive.
The announced result follows as $\theta\leq (8k_{*}K_{W})^{-1}$.
\end{proof}

\begin{proof}[Proof of Lemma \ref{GapBlockLemma}]
Let $\Lambda\geq \Lambda_{0}$. We shall separate the proof in three exhaustive cases.

\noindent
\textbf{Case 1 :} If $\nu_{k+n_{\Lambda}}\leq \Lambda$. Then from Lemma \ref{BelowTheCut} we have $\nu_{k+n_{\Lambda}}-\nu_{k}=\mu_{k+n_{\Lambda}}-\mu_{k}$ and the result follows from Lemma \ref{GapMu}.

\noindent
\textbf{Case 2 :} If $\nu_{k}\leq\Lambda<\nu_{k+n_{\Lambda}}$. Then from Lemma \ref{BelowTheCut} $\nu_{k}=\mu_{k}$ and
from Lemma \ref{BoundForMu}, $\nu_{k+n_{\Lambda}}-\nu_{k} \geq  \mu_{k+n_{\Lambda}}-\mu_{k}$, and thus the result follows from Lemma \ref{GapMu}.

\noindent
\textbf{Case 3 :} If $\Lambda<\nu_{k}<\nu_{k+n_{\Lambda}}$. Then, by definition and Lemma \ref{BelowTheCut}
\[
\nu_{k+n_{\Lambda}} - \nu_{k} = \kappa^{2}\Lambda^{1-2a}((k+n_{\Lambda})^{2}-k^{2}) =  \kappa^{2}\Lambda^{1-2a}((2kn_{\Lambda})^{2}+n_{\Lambda}^{2}) ,
\]
which proves the result.
\end{proof}

\subsubsection{Weak gap property}

\begin{Lemma}\label{WeakGapLemma}
There exists $c=c(c_w,a,\kappa)>0$ such that for any $\Lambda\geq 1$,
\[
\nu_{k+1}-\nu_{k}\geq \exp(-c\Lambda^\frac{1}{2}).
\]
\end{Lemma}
\begin{proof} We shall separate three exhaustive cases.\\

\noindent
\textbf{Case 1 :} If $\mu_{k+1}\leq \Lambda$. Then the result follows from Lemma \ref{BelowTheCut} the assumptions on $(\mu_{k})_{k\in\/N^{*}}$ :
\[
\nu_{k+1}-\nu_{k} = \mu_{k+1}-\mu_{k}\geq\exp(-c_{w}\mu_{k}^{\frac{1}{2}}) \geq\exp(-c_{w}\Lambda^{\frac{1}{2}}) . 
\]

\noindent
\textbf{Case 2 :} If $\mu_{k}> \Lambda$. Then, from Lemma \ref{BelowTheCut}, we have 
\[
\nu_{k+1}-\nu_{k} = \kappa^{2}\Lambda^{1-2a}\geq\exp(-c_0\Lambda^{\frac{1}{2}}), 
\]
with $c_0=  \sup_{x\geq 1}\{x^{-\frac{1}{2}}\ln(\kappa^{-2}x^{2a-1})\}.$ 
\noindent
\textbf{Case 3 :} If $\mu_{k}\leq \Lambda<\nu_{k+1}$. In that case, $k=k^{*}_{\Lambda}$, and
\[
\nu_{k_{\Lambda}^{*}+1}-\nu_{k_{\Lambda}^{*}} = \kappa^{2}\Lambda^{2a-1}(k^{*}_{\Lambda}+1)^{2}-\mu_{k^{*}_{\Lambda}}.
\]
From Lemma \ref{BoundForMu},
\[
\nu_{k_{\Lambda}^{*}+1}-\nu_{k_{\Lambda}^{*}} \geq \mu_{k_{\Lambda}^{*}+1}-\mu_{k_{\Lambda}^{*}} \geq \exp (- c_{w}\mu^\frac{1}{2}_{k^{*}_{\Lambda}})\geq \exp (- c_{w}\Lambda^\frac{1}{2}).
\]
The conclusion follows by setting $c=\max(c_0,c_w).$
\end{proof}

\subsubsection{Proof of Proposition \ref{BiorthoAuxilliary}}\label{ProofOfProposition22}\begin{proof}
    
Let $\Lambda\geq\Lambda_0$. We apply Proposition \ref{GBO} to the sequence $(v_{k})_{k\in\N^{*}}=(\nu_{k})_{k\in\N^{*}}$. Indeed, from Lemma \ref{BoundsOnCounting}, assumption \eqref{WeylGBO} holds with $\alpha=2 c_1\mu_1^\frac{1}{2}$ and $\mathfrak p_{1} = c_{1}$ and $\mathfrak p_{2}=c_{2}\Lambda^{\frac{2a-1}{2}}$. Then we can choose $\mathfrak p_{0}=\min(c_{1},c_{2})>0$, and then from Lemma \ref{EstimatesOnConstants},
\begin{equation}\label{estimatesonrelevantquantities}
(\mathfrak p_{1}\gamma_{\Lambda})^{-2} \leq (c_{1}\theta)^{-2} \Lambda^{\frac{1}{2}-\varepsilon},\quad \mathfrak p_{2} \leq c_{2}\Lambda^{2a-1}\leq c_{2}\Lambda^{\frac{1}{2}-\varepsilon}.
\end{equation}
Moreover, assumption \eqref{GapsGBO} is fulfilled from Lemmata \ref{GapBlockLemma} and \ref{WeakGapLemma} with $\mathfrak n = n_{\Lambda}$, $\mathfrak g=\gamma_{\Lambda}$ and $a_{k}=\exp(c\Lambda^\frac{1}{2})$. From Proposition \ref{GBO}, for any $T\in(0,1)$, there exists $(\sigma_{k})_{k\in\N^{*}}\subset L^{2}(0,T)$ and a constant $C>0$ depending only on $\mu_{1}$, $\min(c_{1},c_{2})$, such that for any $k\in\N^{*}$,
\begin{multline*}
\|\sigma_{k}\|_{L^{2}(0,T)}  \leq C a_{k} ^{n_{\Lambda}}\left( \frac{16+4\gamma_{\Lambda} c_{2}\Lambda^{2a-1}}{\gamma_{\Lambda}^{2} c_{1}^{4}}\right)^{2n_{\Lambda} -2} \\ 
\times \exp\left(C\big(-v_{k}T+v_{k}^{\frac{1}{2}} \left( 1+n_{\Lambda} + 16 (c_{1}\gamma_{\Lambda} )^{-2}+c_{2}\Lambda^{\frac{2a-1}{2}} \right)+\frac{(1+c_{2}\Lambda^{\frac{2a-1}{2}})^{2}}{T}\big)\right).
\end{multline*}
From \eqref{estimatesonrelevantquantities}, we deduce the existence of $C_{0}>0$ (independent on $\Lambda$ and $T\in(0,1)$) such that
\begin{equation*}
\|\sigma_{k}\|_{L^{2}(0,T)}  \leq C_{0} \exp\left( C_{0}\left(\Lambda^{1-\varepsilon}+\nu_{k}^{\frac{1}{2}}\Lambda^{\frac{1}{2}-\varepsilon}+\frac{\Lambda^{1-2\varepsilon}}{T}\right)\right),\quad k\in\N^{*},\; \Lambda\geq \Lambda_{0}.
\end{equation*}
The result follows from the Young inequality.
\end{proof}

\section{The Lebeau-Robbiano method}\label{sec_abstr}
Let $H$ be a real separable Hilbert space endowed with a scalar product  $\langle \cdot,\cdot\rangle_H$ associated to the norm $\|\cdot\|_H$. Let $A:D(A)\subset H\to H$ be a densely defined linear operator such that

\begin{itemize}
\item $A$ is self-adjoint;
\item  $\langle Ax,x\rangle_H \geq 0$ for every $x\in D(A)$;
\item there exists $\mu>0$ such that the operator $(\mu I+A)^{-1}:H\to H$ is compact.
\end{itemize}
In this framework, the spectrum of $A$ is purely discrete and consists of a sequence of ordered non-negative real numbers $(\mu_k)_{k\in\N^*}$. We denote by $(\varphi_k)_{k\in\N^*}$ the associated eigenfunctions, which form a Hilbert basis of $H$. The operator $-A$ generates a strongly continuous semigroup denoted by $e^{-tA}$. 

We now introduce the linear controlled problem
\begin{equation}\label{lin}
\begin{cases}
\xi'(t)+A\xi(t)+Bv(t)=0, &t\in(0,T),\\
\xi(0)=\xi_0\in H,
\end{cases}
\end{equation}
where $v\in L^{2}(0,T)$ and $B\in H$.
We denote by $\xi$ the (mild) solution to the linear control problem \eqref{lin} corresponding to the initial condition $\xi_0$ and control $v$ given by
\[
\xi(t) = e^{-tA}\xi_{0}+\int_{0}^{t} e^{-(t-s)A} Bv(s)ds.
\]
We now introduce the partial null-controllability definition we shall use in what follows. 
\begin{Definition}
Let $E$ be a closed subspace of $H$. We say that the problem \eqref{lin} is null-controllable in $E$ and in a time $T>0$, if there exists a constant $K_{E}(T)>0$ such that, for any $\xi_0\in H$, there exists a control $v\in L^2(0,T)$ such that
$$\Pi_{E}\xi(T)=0\quad\text{ and }\quad\|v\|_{L^2(0,T)}\leq K_{E}(T)\|\xi_0\|_H,$$
where $\Pi_{E}$ denotes the orthogonal projection on $E$ in $H$. The smallest constant $K_{E}(T)$ is called the control cost, and it is defined as
\begin{equation*}
K_{E}(T):=\sup_{\| \xi_0\|_H=1}\inf_{v\in L^{2}(0,T)} \left\{\|v\|_{L^2(0,T)}\,:\,\Pi_{E}\xi(T)=0\right\}.
\end{equation*}
When $E=H$, we simply say that the problem \eqref{lin} is null-controllable and we write $K_{E}(T)=K(T)$.
\end{Definition}

\subsection{Null-controllability in filtered spaces using the moment method}\label{finite}
 Let $T>0$ and $\Lambda>0$. We denote by $$ J_{\Lambda}:=\{k\in\N^*,\;\mu_{k}\leq \Lambda\},\quad\quad\quad E_{\Lambda}:=\Span\{\varphi_k,\; k\in J_{\Lambda}\}.$$
We introduce a classical result of null-controllability on finite-dimensional spaces. 

\begin{Proposition}\label{from_biortho_to_control} Let $\Lambda>0$ and $j\in\N^*$. Assume that $B\in H$ satisfies the following identities: 
$$\langle B,\varphi_{k}\rangle_{H}\neq 0,\ \ \ \forall k\in J_\Lambda.$$
If there exists $(\sigma_{\Lambda,j})_{j\in J_{\Lambda}}\subset L^{2}(0,T)$ such that $$\int_{0}^{T} e^{\mu_{k}t}\sigma_{\Lambda,j}(t)dt=\delta_{j,k},\ \ \ \ \ \ \ \forall j ,k\in J_{\Lambda},$$ where $\delta_{j,k}$ denotes the Kronecker symbol, then the problem \eqref{lin} is null-controllable in $E_{\Lambda}$ in any time $T>0$ with control cost 
\[
K_{E_\Lambda}(T)^2 \leq  (\#J_\Lambda) \frac{\sup_{k\in J_\Lambda} \|\sigma_k\|_{L^2(0,T)}^2}{\inf_{k\in J_\Lambda} |\langle B,\varphi_k\rangle_H|^2}   
\]
\end{Proposition}
\begin{proof}
By the Duhamel formula, the problem \eqref{lin} is null-controllable in $E_{\Lambda}$ and in a time $T>0$, if there exists a function $v_{\Lambda}\in L^{2}(0,T)$ such that, for any $k\in J_{\Lambda}$, the following identities are verified:
\[
\int_{0}^{T}e^{t\mu_{k}} v_{\Lambda}(t) dt = d_{k}, \quad\quad \quad \quad \forall k\in  J_\Lambda ,
\]
where $d_{k}=-\tfrac{\langle \xi_{0},\varphi_{k}\rangle_{H}}{\langle B,\varphi_{k}\rangle_{H}}$. For every $T>0$, we can define $v_{\Lambda}$ as the linear combination
\[
v_{\Lambda} = \sum_{k\in J_{\Lambda}} d_{k}\sigma_{\Lambda,k}.
\]
The estimate of $v_\Lambda$ follows by noticing
\begin{align*}
\| v_{\Lambda} \|^2_{L^2(0,T)}  \leq (\#J_\Lambda) \sum_{k\in J_{\Lambda}} |d_{k}|^2 \|\sigma_{\Lambda,k}\|_{L^2(0,T)}^2  \leq (\#J_\Lambda) \frac{ \sup_{k\in J_\Lambda}\|\sigma_{\Lambda,k} \|_{L^2(0,T)}^2}{\inf_{k\in J_\Lambda}| \langle B,\varphi_{k}\rangle_{H}|^2 }  \| \xi_0 \|^2_H .
\end{align*}
\end{proof}

\subsection{Null-controllability in the full space}\label{lebeau}
In this section, we present the adaptation due to Miller \cite{Miller} (see also \cite{BeauchardPravdaStarov}) of the well-known Lebeau-Robbiano method to prove the observability of the problem \eqref{lin}. First, we need to introduce the definition of partial observability, which is an adaptation of \cite[Proposition 7.7]{LRLR} to our case. 
\begin{Proposition}\label{PartialObs}Let $E$ be a closed subspace of $H$. The null-controllability of the problem \eqref{lin} in $E$ at time $T>0$ with control cost $K_{E}(T)$ is equivalent to the following observability estimate : 
\[
\| e^{-TA}\Pi_{E} \zeta \|_{H} \leq K_{E}(T)\| \langle B,e^{-tA}\Pi_{E} \zeta\rangle_{H} \|_{L^{2}(0,T)}, \quad \forall \zeta\in H.
\]
\end{Proposition}
\begin{proof}
Consider the bounded control operator, with $U=\R$,
\[
\begin{array}{rccl}
M : &U& \rightarrow &H \\ &v &\mapsto & Bv,
\end{array}
\]
whose adjoint is given by 
\[
\begin{array}{rccl}
M^{*} : &H& \rightarrow &U \\ & z &\mapsto & \langle B,z \rangle_{H} .
\end{array}
\]
The proof follows by a direct application of \cite[Proposition 7.7]{LRLR}.
\end{proof}

We are finally ready to show how to use the Lebeau-Robbiano-Miller method to prove the null-controllability of the problem \eqref{lin} when partial observability inequalities are verified in suitable finite-dimensional subspaces of $H$ and the control cost satisfies specific inequalities.

\begin{Theorem}\label{LRM}
Assume that there exist $C_{0}>0$ and $\alpha\in(0,1)$ such that for any $T\in(0,1)$ and for any $\Lambda>0$, the problem \eqref{lin} is null-controllable in $E_{\Lambda}$ at time $T>0$ with control cost 
\[
K_{E_{\Lambda}}(T)\leq C_{0}\exp\left(C_{0}\big(\Lambda^\alpha+\frac{1}{T^\frac{\alpha}{1-\alpha}}\big)\right).
\]
Then there exists $C_{1}>0$ such that for any $T\in (0,1)$, the problem \eqref{lin} is null-controllable with control cost 
\[K(T)\leq C_{1}\exp \left(\frac{C_{1}}{T^\frac{\alpha}{1-\alpha}}\right).\]
\end{Theorem}
\begin{proof}
Let $\tau>0$ and $0<T_{1}<T_{2}<1$ be such that $\tau=T_{2}-T_{1}$. Let $\delta\in(0,1/2)$ to be fixed below. We have, for any $\zeta\in H$, and any $\Lambda>0$,
\begin{equation*}
\|e^{-T_{2}A}\zeta \|_{H}^{2} = \|e^{-T_{2}A}\Pi_{E_{\Lambda}}\zeta \|_{H}^{2} + \|e^{-T_{2}H}(I-\Pi_{E_{\Lambda}})\zeta \|_{H}^{2} .
\end{equation*}
From Proposition \ref{PartialObs}, it follows
\begin{equation*}
\|e^{-T_{2}A}\zeta \|_{H}^{2} \leq  K_{E_{\Lambda}}(\delta \tau)^{2}\| \langle B,e^{-2tA}\Pi_{E_{\Lambda}} \zeta\rangle_{H} \|^{2}_{L^{2}(T_{2}-\delta \tau,T_{2})} + e^{-\Lambda \tau} \|e^{-2T_{1}A}(I-\Pi_{E_{\Lambda}})\zeta \|_{H}^{2} .
\end{equation*}
Yet, using the Cauchy-Schwarz inequality
\begin{align*}
\|& \langle B,e^{-tA}\Pi_{E_{\Lambda}} \zeta\rangle_{H} \|_{L^{2}(T_{2}-\delta \tau,T_{2})}^{2} \\
 \leq &\;2 \| \langle B,e^{-tA}\zeta\rangle_{H} \|_{L^{2}(T_{1},T_{2})}^{2}+ 2 \| \langle B,e^{-tA}(I-\Pi_{E_{\Lambda}}) \zeta\rangle_{H} \|_{L^{2}(T_{2}-\delta \tau,T_{2})}^{2} \\
   \leq &\; 2 \| \langle B,e^{-tA}\zeta\rangle_{H} \|_{L^{2}(T_{1},T_{2})}^{2}+2 e^{-2(1-\delta)\Lambda\tau} \|  B\|^{2}_{H}\| e^{-T_{1}A}\zeta \|_{H}^{2}.
\end{align*}
The above estimates yield the existence of $C>0$ such that, for any $0<\tau$, $0<T_{1}<T_{2}<1$ such that $\tau=T_{2}-T_{1}$, $\delta\in(0,\tfrac{1}{2})$, $\Lambda>0$, and any $\zeta\in H$,
\begin{multline*}
\|e^{-T_{2}A}\zeta \|_{H}^{2} \leq  C\exp \left(C(\Lambda^\alpha+(\delta \tau)^\frac{\alpha}{\alpha-1})\right) \| \langle B,e^{-tA}\zeta\rangle_{H} \|_{L^{2}(T_{1},T_{2})}^{2}\\
+\Big(C\exp\left(C(\Lambda^\alpha+(\delta \tau)^\frac{\alpha}{\alpha-1})-2\Lambda(1-\delta)\tau\right) + \exp(-2\Lambda \tau) \Big) \|e^{-T_{1}A}\zeta \|_{H}^{2} .
\end{multline*}
Now, we choose $\Lambda=(\delta\tau)^{\frac{1}{\alpha-1}}$ and there exists $\tilde C>0$ independent on $\delta$, $\tau$, and $\zeta$ such that
\begin{multline*}
\|e^{-T_{2}A}\zeta \|_{H}^{2} \leq  \tilde C\exp\left(\tilde C(\delta \tau)^\frac{\alpha}{\alpha-1}\right) \| \langle B ,e^{-tA}\zeta\rangle_{H} \|_{L^{2}(T_{1},T_{2})}^{2}\\+
\Big( \tilde C\exp\left(\tilde C(\delta \tau)^\frac{\alpha}{\alpha-1}-2(1-\delta)\delta^{\frac{1}{\alpha-1}}\tau^\frac{\alpha}{\alpha-1}\right)+ \exp\left(-2\delta^{\frac{1}{\alpha-1}}\tau^\frac{\alpha}{\alpha-1}\right)\Big) \|e^{-T_{1}A}\zeta \|_{H}^{2} .
\end{multline*}
Let $f_{\delta} (\tau) = \tilde C^{-1}\exp(-\tilde C(\delta\tau)^\frac{\alpha}{\alpha-1})$. We fix $\delta\in(0,1)$ sufficiently small (independently on $\tau,T_1,T_2$) so that the following approximate observability holds : there exists $\tilde C>0$ such that for any $0<T_{1}<T_{2}<1$, we have
\[
f_{\delta}(\tau) \|e^{-T_{2}A}\zeta \|_{H}^{2}  \leq  \| \langle B ,e^{-tA}\zeta\rangle_{H} \|_{L^{2}(T_{1} ,T_{2})}^{2} + f_{\delta}\Big(\frac{\tau}{2}\Big) \|e^{-T_{1}A}\zeta \|_{H}^{2}  .
\]
For any $T\in (0,1)$ and $k\in\N$, we choose $T_{2}=\tfrac{T}{2^{k}}$ and $T_{1}=\tfrac{T}{2^{k+1}}$ with $k\in\N$, and we obtain 
\begin{align*}
f_{\delta}(\tfrac{T}{2^{k}}) \|e^{-T/2^{k}A}\zeta \|_{H}^{2}  
\leq  \| \langle B,e^{-tA}\zeta\rangle_{H} \|_{L^{2}\big(\tfrac{T}{2^{k+1}} ,\tfrac{T}{2^{k}}\big)}^{2} + f_{\delta}(\tfrac{T}{2^{k+1}}) \|e^{-\tfrac{T}{2^{k+1}}A}\zeta \|_{H}^{2}.  
\end{align*}
Summing over $k\in\N$ yields the following observability inequality 
\[
\|e^{-TA}\zeta \|_{H}^{2}  \leq f_{\delta}^{-1}(T)  \| \langle B,e^{-tA}\zeta\rangle_{H} \|_{L^{2}(0 ,T)}^{2}.  
\]
Finally, it is sufficient to apply Proposition \ref{PartialObs} to end the proof.
\end{proof}

\section{Proof of the main results}
\subsection{Proof of Theorem \ref{NullControlTheorem}}\label{ProofOfTheSecondTheorem}

From Propositions \ref{BiorthoAuxilliary} and \ref{from_biortho_to_control}, we deduce that for any $T>0$ and any $\Lambda>\Lambda_0$ (defined in \eqref{DefsOfGapConstants}), \eqref{lin2} is null-controllable in $E_\Lambda$ with a control cost of the form
\[
K_{E_\Lambda}(T)^2\leq \frac{C}{\inf_{k\in J_\Lambda}\{ |\langle \mathcal B,\phi_k\rangle_{\mathcal H}|^2\}} \exp \left(C(\Lambda^{1-\varepsilon}+\frac{1}{T^{\alpha_{\varepsilon}}})\right),
\]
with $C>0.$ Note that, up to changing the constant $C>0$, \eqref{lin2} is null-controllable in $E_\Lambda$ for any $\Lambda>0$ with the same form of control cost. Moreover, from \eqref{AssumptionOnB} and \eqref{WeylLaw1Term}, 
\[
\inf_{k\in J_\Lambda}\{ |\langle \mathcal B,\phi_k\rangle_{\mathcal H}|^2\} \geq \exp\left(-C_{\mathcal B} \Lambda^{1-\delta }\right) 
\]
Then it is sufficient to apply Theorem \ref{LRM} to obtain the announced result.

\subsection{Proof of Theorem \ref{SolvabilityOfTheMoment}}\label{ProofOfTheMainResult}
We place ourselves in the framework of Section \ref{NullControlLinearSection}. Let $x=(x_k)_{k\in\N^*}\in \ell^2$ and $z=(\langle\mathcal B, \phi_k\rangle_{\mathcal H})_{k\in\N^*}$ be such that $zx\in\ell^2$. We define
\[
\xi_0=-\sum_{k\in\N^*} z_kx_k \phi_k\in \mathcal H.
\]
From Theorem \ref{NullControlTheorem}, for any $T\in(0,1)$, there exists $v\in L^2(0,T)$ satisfying 
\[
\|v\|^2_{L^2(0,T)} \leq C\exp\left( \frac{C}{T^{\alpha_{\varepsilon,\delta}}} \right)\|\xi_0\|^2 =  C\exp\left( \frac{C}{T^{\alpha_{\varepsilon,\delta}}} \right)\|zx\|_{\ell^2}^2,
\]
and such that the mild solution of \ref{lin2} satisfies for any $T>0$,
\[
 -e^{-T\mathcal A}\xi_0 = \int_0^T e^{-(T-t)\mathcal A}\mathcal Bv(t)dt,
\]
which is equivalent to
\[
\int_0^Te^{t\mu_k} v(t)dt=-\frac{\langle \xi_0,\phi_k\rangle_{ \mathcal H}}{\langle \mathcal B,\phi_k\rangle_{\mathcal H}} =x_k, \quad \forall k\in\N^*.
\]
This ends the proof.

\section{Local exact controllability to the eigensolutions}
\label{Section_Control_Simple_Spectrum}

\subsection{Proof of Theorem \ref{BilinearControlTheo}}
 The aim of this section is to prove the exact controllability to eigensolutions of the bilinear parabolic equation \eqref{main}. Before moving on with the proof of Theorem \ref{BilinearControlTheo}, we provide the well-posedness result of \eqref{main}. It is a classical result that we rephrase in the next proposition (see for instance \cite{bms}).

\begin{Proposition}\label{well-posedness}
Let $T>0$. For any $\psi_0\in \mathcal H$, there exists a unique mild solution $\psi\in C^0([0,T],\mathcal H)$ to \eqref{main} such that
\begin{equation*}
    \psi(t)=e^{-t \mathcal A}\psi_0-\int_0^t e^{-(t-s)\mathcal A}v(s)\mathcal B \psi(s)ds,\quad\forall t\in[0,T].
\end{equation*}
Moreover, there exists a constant $C(T)>0$ such that $\sup_{t\in[0,T]}\|\psi(t)\|_{\mathcal H}\leq C(T)\|\psi_0\|_{\mathcal H}.$ 
\end{Proposition}
 We introduce the following linearized problem in $\mathcal H$
\begin{equation}\label{lin_gen}
\begin{cases}
\xi'(t ) +\mathcal A \xi(t ) + v(t) \mathfrak B \phi_{j}  =0, & t\in(0,T) ,\\
\xi(0 )=\xi_0\in\mathcal H.
\end{cases}
\end{equation}
The proof of Theorem \ref{BilinearControlTheo} is a direct consequence of the following result.

\begin{Proposition}\label{prop:local-nul}
If the problem \eqref{lin_gen} is null-controllable in any $T>0$  and there exist some constants $\nu,T_0>0, \gamma\geq1$ such that 
\begin{equation}\label{bound-control-cost}
K(\tau)\leq e^{\nu/\tau^\gamma},\quad\forall\,0<\tau\leq T_0.
\end{equation}
then, System \eqref{main} is small-time locally exactly controllable to eigensolutions.\end{Proposition}

Proposition \ref{prop:local-nul} is a standard result that follows from the arguments of \cite[Theorem 1.1]{acue}. The derivation in \cite[Section 3]{acue} handles the case $\gamma=1$, but the proof extends naturally to any $\gamma \geq 1$. Such results are quite standard in the literature (see also \cite{BeauchardMarbach,LTT}), and the proof is omitted for brevity.

\begin{proof}[Proof of Theorem \ref{BilinearControlTheo}]
Let us consider \eqref{main}.
We observe that, under the hypotheses of Theorem \ref{BilinearControlTheo}, the problem \eqref{lin_gen} is null-controllable thanks to Theorem \ref{NullControlTheorem} (or Corollary \ref{coro:control}). Moreover, the corresponding control cost satisfies an upper bound of the form \eqref{bound-control-cost}. Finally, the exact controllability to the eigensolutions of \eqref{main} stated by Theorem \ref{BilinearControlTheo}, follows from Proposition \ref{prop:local-nul}.
\end{proof}

\subsection{Bilinear controllability on a rectangle}

Let us now discuss the specific case $\mathcal H = L^2(\Omega)$, where
$$\Omega=(0,a)\times (0,b),$$
with $a,b>0$ such that $a^2/b^2$ is an algebraic irrational number. In this case, we consider $\mathcal A$ to be the Bi-Laplace operator $\Delta^2$ with hinged boundary conditions, with domain
$$D(\Delta^2):=\{\psi \in H^4(\Omega): \; \psi,\Delta\psi \in H^1_0(\Omega)\}.$$
Let $Q^1\in L^\infty(0,a)$, $Q^2\in L^\infty(0,b)$ such that, for every $\ell\in\N^*$, there exist $C>0$ and $p_1,p_2> 0$ such that for any $k\in\N^*$
\begin{equation}\begin{split}\label{AssumptionsOnBrect} \big| \int_0^a\sin\big(\tfrac{k}{a}\pi x\big) Q_1(x)\sin\big(\tfrac{\ell }{a}\pi x\big) dx\big|\geq \frac{C}{k^{p_1}},\\
\big| \int_0^b \sin\big(\tfrac{k}{b}\pi y\big) Q^2(y)\sin\big(\tfrac{\ell }{b}\pi y\big)dy\big| \geq \frac{C}{k^{p_2}}.\end{split}\end{equation}
 We define $Q(x,y)=Q^1(x)Q^2(y)\in L^\infty(\Omega)$, and $\mathfrak B$ as the multiplication operator 
$$
\begin{array}{rcl} 
\mathfrak B: L^2(\Omega) & \longrightarrow & L^2(\Omega) \\
\psi  & \longmapsto & Q\psi ,
\end{array}$$ 
In this case, System \eqref{main} reads 
\begin{equation}\label{main-spec}
\begin{cases}
\partial_t \psi +\Delta^2\psi+v Q \psi=0, & \text{ in }(0,T)\times\Omega,\\
\psi_{|_{\partial\Omega}}=\Delta \psi_{|_{\partial\Omega}} =0,& \text{ in }(0,T)\times\partial\Omega\\
\psi_{|_{t=0}}=\psi_0\in L^2(\Omega),& \text{ in } \Omega,
\end{cases}
\end{equation}
with $v\in L^2(0,T)$. In this framework, the ordered eigenvalues $\big(\mu_k\big)_{k\in\N^*}$ are defined by two sequences of numbers $(\ell_k)_{k\in\N^*},(m_k)_{k\in\N^*}\subset\N^*$ such that
\begin{align}\label{l_k} \mu_k =\pi^4 \Big(\frac{\ell_k^2}{a^2}+\frac{m_k^2}{b^2}\Big)^2,\end{align}
with corresponding eigenfunctions $\big(\phi_k\big)_{k\in\N^*}$ given by
$$\phi_k(x,y)=\frac{2}{\sqrt{a\, b}} \sin\Big(\frac{\ell_k}{a}\pi x\Big)\sin\Big(\frac{m_k}{b} \pi y\Big) .$$
The suitable choice of $a$ and $b$ so that $a^2/b^2$ is an algebraic irrational number allows to ensure the exact controllability of Theorem \ref{BilinearControlTheo} for \eqref{main-spec}.

\begin{Remark}\label{ex1_rec}
In the case where $\Omega=(0,1)\times\big(0,\sqrt[3]2\big)$, and $Q^1(x)= x^2$, $Q^2(y)=y^2,$ our hypotheses are verified. Indeed, the number $2^{-\frac{2}{3}}$ is algebraic and irrational. In addition, for every $\ell\in\N^*$ and $k\in\N^*\setminus\{\ell\}$, we have
\begin{align*}\big | \int_0^1{t^2} \sin ( {k}  \pi t )   \sin ( {\ell}  \pi t )  dt \big| 
    &= \frac{1}{2}\big| \int_0^{2^\frac{1}{3}} t^2 \sin\Big(2^{-\frac{1}{3}} k\pi t\Big) \sin\Big(2^{-\frac{1}{3}}\ell \pi t\Big)dt \big| \\
    &
    =\frac{4  k\, \ell}{(k^2-m^2)^2 \pi^2} \geq \frac{C_\ell}{k^3}.\end{align*}
When $k=\ell$, we compute
\begin{align*} \big|  {t^2} \int_0^1 \sin^2 ( {\ell}  \pi  ) dt \big|
    = \frac{1}{2}\big|   \int_0^{2^\frac{1}{3}} t^2 \sin^2\Big(\frac{\ell}{2^\frac{1}{3}} \pi t\Big) dt\big|
    =\frac{ \left| 2 \ell^2 \pi^2 - 3 \right| }{12 \ell^2 \pi^2}> 0.\end{align*}
 Similar results can be obtained for other polynomial controls $Q^1$ and $Q^2$ of different orders. 
\end{Remark}

\begin{Theorem}\label{th:local-rec}
System \eqref{main-spec} is small-time locally exactly controllable to eigensolutions.
\end{Theorem}

\begin{proof}
The result is a direct application of Theorem \ref{BilinearControlTheo} to this setting. Indeed, the eigenvalues of the Dirichlet Bi-Laplace operator on a rectangle satisfy \eqref{DefOfMu}, \eqref{WeylLaw}, and \eqref{AssumptionOnAB} with $a=\tfrac{1}{2}$ and $b=\tfrac{1}{4}$ (similarly to the square case described in Section \ref{ExampleSectionMoment}).
The validity of \eqref{Weakgap}, is ensured as follows. First, Roth's Theorem (see \cite{Roth}) implies that, when $z$ is an algebraic irrational number, we have that, for every $\epsilon>0$, there exists $C_\epsilon>0$ small enough such that
$$\Big|z-\frac{n}{m}\Big|\geq\frac{C_\epsilon}{m^{2+\epsilon}},\quad\forall m,n\in\Z^*.$$
Recall \eqref{l_k}. Now, $a^2/b^2$ is an algebraic irrational number and for any $\epsilon>0$, there exists $C_\epsilon>0$ such that
\begin{equation*}\begin{split}
|\mu_{k+1}-\mu_k|&=\Big|\Big(\frac{\ell_{k+1}^2\pi^2}{a^{2}}+\frac{m_{k+1}^2\pi^2}{b^{2}}\Big)^2-\Big(\frac{\ell_k^2\pi^2}{a^{2}}+\frac{m_k^2\pi^2}{b^{2}}\Big)^2\Big|\\
&=\pi^2\Big|\Big(\frac{\ell_{k+1}^2}{a^{2}}+\frac{m_{k+1}^2}{b^{2}}\Big)-\Big(\frac{\ell_k^2}{a^{2}}+\frac{m_k^2}{b^{2}}\Big)\Big|\times \big|\sqrt{\mu_{k+1}}+\sqrt{\mu_k}\big|\\
&=\frac{ \pi^2|m_k^2-m_{k+1}^2|}{a^{2}}\Big|\frac{\ell_{k+1}^2-\ell_k^2}{m_k^2-m_{k+1}^2} -\frac{a^{2}}{b^2}\Big|\times \big|\sqrt{\mu_{k+1}}+\sqrt{\mu_k}\big|\\
&\geq 2\sqrt{\mu_{k+1}} \frac{\pi^2 |m_k^2-m_{k+1}^2|C_\epsilon}{a^2|{m_k^2-m_{k+1}^2}|^{2+\epsilon}}\geq 2\sqrt{\mu_{k+1}}    \frac{\pi^2 C_\epsilon}{a^2 \max(m_k,m_{k+1})^{2+2\epsilon}}\\
&\geq 2\mu_{k+1}^\frac{1}{2} \frac{\pi^{4+2\epsilon} C_\epsilon}{a^2 b^{2+2\epsilon} \mu_{k+1}^\frac{1+\epsilon}{2}}= \frac{2\pi^{4+2\epsilon} C_\epsilon}{a^2 b^{2+2\epsilon} }\frac{1}{\mu_{k+1}^\frac{\epsilon}{2}}.
\end{split}\end{equation*}
The Weyl asymptotics ensures the validity of \eqref{Weakgap}. Finally, we can apply Theorem \ref{BilinearControlTheo}, since \eqref{AssumptionOnFrakB} is ensured by the identities \eqref{AssumptionsOnBrect}.

\end{proof}

\appendix
\section{The case of a spectral sequence with weaker accumulation}\label{WeakerAccumulationSection}
 Here, we define for $a\in(0,\tfrac{1}{2}),$
    \begin{equation}\label{DefOfBeta}
        \beta_{a,\delta}=\max\left(\frac{4}{a(1-2a)},\frac{1-\delta}{\delta} \right)
    \end{equation}
\begin{Corollary}
Let $\mu=(\mu_{k})_{k\in\N^{*}}$ be an increasing sequence of positive numbers satisfying \eqref{DefOfMu}, \eqref{WeylLaw} with $0<b<a<\tfrac{1}{2}$ and \eqref{Weakgap}. Let $z=(z_k)_{k\in\N^*}\in \ell^2$ such that there exists $C_z>0$ such that for any $k\in\N^*$,
\[|z_k|\geq \exp \left( -C_z\mu_k^{1-\delta} \right) .
\]Then, there exists $C:=C(a,b,\varepsilon,\delta,z,\mu,c_{W,1},c_{W,2},c_w)>0$ such that for any $T\in(0,1)$ and for any $x=(x_{k})_{k\in\N^{*}}\subset \R$ satisfying $zx \in\ell^2$, there exists $u\in L^{2}(0,T)$ such that
\[
x_{k}=\int_{0}^{T}e^{t\mu_{k}}u(t) dt,
\]
with the following estimate
\[
\|u\|_{L^{2}(0,T)} \leq C\exp\left( \frac{C}{T^{\beta_{a,\delta}}} \right) \| zx \|_{\ell^{2}}.
\] 
\end{Corollary}
\begin{proof} Let $\mu=(\mu_{k})_{k\in\N^{*}}$ satisfy \eqref{DefOfMu}, \eqref{WeylLaw} with $0<b<a<\frac{1}{2}$ and \eqref{Weakgap}. We then consider an auxilliary sequence $\lambda_{k,\ell}$ defined by 
    \[
\lambda_{k,\ell}=\mu_k+\ell \frac{\mu_{k+1}-\mu_k}{k^\frac{1-2a}{2a}}, \quad k\in\N^*,\; \ell\in \{0,\dots, \lfloor k^{\frac{1-2a}{2a}} \rfloor-1 \}.
    \]
Let $\Gamma\geq \mu_1$. Let $k:=k(\Gamma)\in\N^*$ such that $\mu_k\leq \Gamma < \mu_{k+1}$. Its counting function satisfies 
    \begin{align*}
           \mathcal N_{\lambda}(\Gamma) & =    \mathcal N_{\lambda}(\mu_{k}) + \#\{\ell\in \{0,\dots, \lfloor k^{\frac{1-2a}{2a}} \rfloor-1 \},\; \mu_k< \lambda_{k,\ell}\leq \Gamma\}.
    \end{align*}
     Yet, by definition
    \[
    \#\{\ell\in \{0,\dots, \lfloor k^{\frac{1-2a}{2a}} \rfloor-1 \},\; \mu_k< \lambda_{k,\ell}\leq \Gamma\} \leq k^{\frac{1-2a}{2a}} \lesssim \Gamma^{\frac{1-2a}{2}}.
    \]

    Moreover,
    \begin{align*}
    \mathcal N_{\lambda}(\mu_{k}) & = \sum_{\mu_j\leq \mu_k} \lfloor j^{\frac{1-2a}{2a}} \rfloor =\sum_{1\leq j\leq k} j^{\frac{1-2a}{2a}} -\sum_{1\leq j\leq k} \left( j^{\frac{1-2a}{2a}} - \lfloor j^{\frac{1-2a}{2a}}\rfloor\right) \\
    & = 2a k^{\frac{1}{2a}}+ \mathcal O (k^\frac{1-2a}{2a})+\mathcal O(k)  = 2a (k+1)^{\frac{1}{2a}}+ \mathcal O (k^\frac{1-2a}{2a})+\mathcal O(k) .
    \end{align*}
    We deduce, for any $\Gamma\geq \mu_1$,
    \[
    \left |\mathcal N_\lambda(\Gamma)-2a \Gamma^\frac{1}{2} \right | \leq C \Gamma^{\max(a,\frac{1-2a}{2})}.
    \]
  On the one hand, there exists $c>0$ such that for any $k\in\N^*$ and any $\ell\in \{0,\dots, \lfloor k^{\frac{1-2a}{2a}} \rfloor-2 \}$,
   \[
   \lambda_{k,\ell+1}-\lambda_{k,\ell} \geq   \frac{\mu_{k+1}-\mu_k}{k^\frac{1-2a}{2a}} \geq\exp\left( - c \mu_{k}^{\frac{1}{2}} \right) \geq\exp\left( - c \lambda_{k,\ell}^{\frac{1}{2}} \right).
   \]
    On the other hand, there exists $c>0$ such that
      \begin{multline*}
   \lambda_{k+1,0}-\lambda_{k, \lfloor k^{\frac{1-2a}{2a}} \rfloor-1} =( \mu_{k+1} - \mu_k) \left(1 + \frac{1-\lfloor k^{\frac{1-2a}{2a}} \rfloor}{k^\frac{1-2a}{2a}} \right) 
   \\ 
   \geq \exp\left(-c \mu_k^\frac{1}{2} \right)
   \geq \exp\Big(-c \lambda_{k, \lfloor k^{\tfrac{1-2a}{2a}}\rfloor-1}^{\frac{1}{2}} \Big) .
   \end{multline*}
    Hence, the sequence $\lambda$ now satisfies \eqref{DefOfMu}, \eqref{WeylLaw}, \eqref{AssumptionOnAB} and \eqref{Weakgap}. Now we introduce two sequences, $\widetilde z_{k,\ell}=z_k$, for any $k\in\N^*$ and any $\ell\in \{0,\dots, \lfloor k^{\frac{1-2a}{2a}} \rfloor-1 \}$, and 
   \[
y_{k,\ell}=\begin{cases}
x_k & \ell=0 \\
0 & \ell\neq 0.
\end{cases}
\]
    We can then apply Theorem \ref{SolvabilityOfTheMoment} to the sequences $\lambda$ and $\widetilde z$, and solve 
   \[
y_{k,\ell}=\int_{0}^{T}e^{t\lambda_{k,\ell}}u(t) dt,
\]
with the estimate
\[
\|u\|_{L^2(0,T)} \leq C\exp\left( \frac{C}{T^{\beta_{a,\delta}}}\right) \|\widetilde z y\|_{\ell^2} = C\exp\left( \frac{C}{T^{\beta_{a,\delta}}}\right) \|zx\|_{\ell^2} .
\]
\end{proof}

\begin{Corollary}\label{coro:control}
Let $\mu=(\mu_{k})_{k\in\N^{*}}$ be an increasing sequence of positive numbers satisfying \eqref{DefOfMu}, \eqref{WeylLaw} with $0<b<a<\tfrac{1}{2}$ and \eqref{Weakgap}, and assume \eqref{AssumptionOnB}. There exists $C:=C(a,b,\varepsilon,\delta,\mathcal B,\mu,c_{W,1},c_{W,2},c_w)>0$ such that for any $\xi_0\in \mathcal H$ and any $T\in(0,1)$, there exists $v\in L^2(0,T)$ such that the unique mild solution of \eqref{lin2} satisfies
\[
\xi(T)=0,
\]
with control cost estimate
    \begin{equation*}
    \| v \|_{L^2(0,T)} \leq C\exp \left(\frac{C}{T^{\beta_{a,\delta}}}\right) \|\xi_0\|_{\mathcal H},
    \end{equation*}
    where $\beta_{a,\delta}$ is defined in \eqref{DefOfBeta}.
\end{Corollary}

\bibliographystyle{abbrv}
\bibliography{refs.bib} 
	\end{document}